\documentclass[12pt]{article}

\usepackage{amsmath}
\usepackage{amssymb}

\newtheorem{Th}{Theorem}[section]

\newtheorem{lemma}[Th]{Lemma}
\newtheorem{corollary}[Th]{Corollary}
\numberwithin{equation}{section}

\newcommand{\finishproof}{\hfill $\Box$ \vspace{3mm}}
\newcommand{\re}{\mathop{\rm Re}}

\begin{document}

\title{Interpolation of nonlinear maps}

\author{T. Kappeler\footnote{Supported in part by the Swiss National Science Foundation.},
A. Savchuk\footnote{ Supported in part by the Russian Fund of Fundamental Research, \# 13-01-00705.},
A. Shkalikov\footnote{  Supported in part by the Russian Fund of Fundamental Research, \# 13-01-00705.},
P. Topalov\footnote{Supported in part by NSF DMS-0901443.}}

\maketitle

\begin{abstract}
\noindent Let $(X_0, X_1)$ and $(Y_0, Y_1)$ be complex Banach couples and assume that $X_1\subseteq X_0$  with norms
satisfying $\|x\|_{X_0} \le c\|x\|_{X_1}$ for some $c > 0$. For any $0<\theta <1$, denote by
$X_\theta = [X_0, X_1]_\theta$ and $Y_\theta = [Y_0, Y_1]_\theta$ the complex interpolation spaces
and by $B(r, X_\theta)$, $0 \le \theta \le 1,$ the open ball of radius $r>0$ in $X_\theta$,
centered at zero.  Then for any analytic map
$\Phi: B(r, X_0) \to Y_0+ Y_1$ such that $\Phi: B(r, X_0)\to Y_0$ and $\Phi: B(c^{-1}r, X_1)\to Y_1$ are continuous and bounded
by constants $M_0$ and $M_1$, respectively, the restriction of $\Phi$ to $B(c^{-\theta}r, X_\theta)$, $0 < \theta < 1,$ is shown
to be a map with values in $Y_\theta$ which is analytic and bounded by $M_0^{1-\theta} M_1^\theta$.

\vspace{0.5cm}

\end{abstract}

\section{Introduction}\label{Introduction}

\noindent Let us first recall some basic notations and definitions of the interpolation theory for Banach spaces.
Following \cite{Tr}, we say that two complex Banach spaces $X_0$, $X_1$  are a complex Banach couple or  Banach couple for short, $(X_0, X_1)$, if they are
both linearly and continuously embedded into a linear complex Hausdorff space $\mathcal X$, i.e., $X_0\subseteq{\mathcal X}$ and
$X_1\subseteq{\mathcal X}$.
The spaces $X_\cap := X_0\cap X_1$ and
$$
X_+ := X_0 + X_1= \{ x\in{\mathcal X}\,|\, x= x_0+x_1, \, x_0\in X_0, \, x_1\in X_1\}
$$
are also Banach spaces when endowed with the norms
$$
\|x\|_\cap := J(1, x) \quad \mbox{ respectively } \quad \|x\|_+ := K(1, x)
$$
where for any $t \ge 0$,
$$
J(t,x) := \max\big\{\|x\|_{X_0}, t\|x\|_{X_1}\big\} \quad \forall \, x \in X_\cap
$$
and for any $x \in X_+,$
\begin{equation}
\label{1}
 K(t,x) := \inf\big\{ \|x_0\|_{X_0} + t\|x_1\|_{X_1}\,\big|\, x=x_0+x_1, \, x_0\in X_0, \, \ x_1\in X_1\big\} .
\end{equation}
Without loss of generality, we will assume in the sequel that $\mathcal X = X_+.$
We say that $(X_0,X_1)$ is a regular Banach couple if  $X_1$ is continuously embedded into $X_0$, i.e.,
$ X_1\subseteq X_0$ and there exists $c > 0$ so that
$$
 \|x\|_{X_0}  \leqslant c \|x\|_{X_1} \quad \forall \, x \in X_1.
$$
There are several methods to construct for any given Banach couple $(X_0, X_1)$ interpolation
spaces $X$, $X_\cap \subseteq X \subseteq X_+$, that satisfy the interpolation property discussed below.
The most familiar interpolation methods are the real and the complex ones.
The real method comprises among others the $K$-method, the $L$-method, the $J$-method,
the mean-methods, and the trace method -- see e.g. \cite{BL}, \cite{Tr}. Up to equivalent norms, these methods all lead to the
same interpolation spaces $X_{\theta, p}$ where $0 < \theta < 1$ and $1 \le p \le \infty$.
In particular,  for any fixed $1\le p < \infty$  and $0< \theta <1$,  the $K$-method defines the interpolation spaces
$X_{\theta, p} \equiv (X_0,X_1)_{\theta, p}$ as follows
\begin{equation}\label{K}
X_{\theta, p} := \Big\{ x\in X_+ \,\Big|\, \|x\|_{\theta, p} =
\Big(\int_0^\infty\big(t^{-\theta} K(t, x)\big)^p \frac {dt}t \Big)^{1/p}<\infty\Big\}
\end{equation}
where $K(t,x)$ is defined by \eqref{1}.  If $p=\infty$, the space $(X_0, X_1)_{\theta,\infty}$
consists of the elements $x \in X_+$ with $\|x\|_{\theta,\infty}= \sup_{0<t<\infty} t^{-\theta} K(t,x) <\infty$.

\noindent Now let us turn to the complex method of interpolation.
For the notion of an analytic map between complex Banach spaces
we refer the reader e.g. to  \cite{Die} or \cite[Appendix A]{PT}.
Following \cite{Tr},  denote by $S$ the vertical strip in the complex plane given by
$$
S :=\{z \in \mathbb C \vert \  0<\re z <1\}
$$
and by $\overline S$ its closure.  For any given Banach couple $(X_0, X_1)$, we then introduce
the complex vector space
$\mathcal H(X_0, X_1)$ of maps $f: {\overline S} \to X_+$ with the following properties:
\begin{itemize}
\item[($\mathcal H 0$)] $f : \overline S \to X_+$ is continuous and bounded;
\item[($\mathcal H 1$)] $f|_{S}: S \to X_+$ is analytic;
\item[($\mathcal H 2$)] for any $t \in \mathbb R$, $f(it) \in X_0$, $f(1+it) \in X_1$, and the maps
$\mathbb R \to X_0, t \mapsto f(it)$ and $\mathbb R \to X_1, t \mapsto f(1 + it)$ are bounded and continuous.
\footnote{The continuity of these maps is not assumed in \cite[Appendix to IX.4]{RS}.
Using the definition of the space $\mathcal H$ from \cite[Appendix to IX.4]{RS} one gets a slightly stronger version of
Theorem \ref{th:main} in the sense that the assumption of the continuity of the maps
\eqref{est0} and \eqref{est1}  can be dropped.}
\end{itemize}
Then
$$
\|f\|_{\mathcal H} := \max \big\{\sup_{t\in \mathbb R}\|f(it)\|_{X_0}, \  \sup_{t\in \mathbb R}\|f(1+it)\|_{X_1}\big\}
$$
defines a norm on $\mathcal H(X_0, X_1)$.
For any $0 < \theta < 1$, the complex interpolation space $X_\theta := [X_0, X_1]_\theta $ is the space
$$
X_\theta := \{  x\in X_+\,|\,  \exists \, f(z)\in \mathcal H(X_0, X_1)\ \text{with}\ f(\theta) = x\},
$$
endowed with the norm
\begin{equation}\label{eq:theta-norm}
\|x\|_{X_\theta} := \inf_{f\in \mathcal{H}} \{\|f\|_{\mathcal H}\,|\,f(\theta) = x\}.
\end{equation}
It is well known that the  spaces $X_\theta$ (and $ X_{\theta, p}$) have the {\em interpolation property}:
given any two Banach couples $(X_0,X_1)$ and $(Y_0, Y_1)$ and any bounded linear operator $T : X_+ \to Y_+$ such
that for some $M_0 > 0$, $M_1 > 0,$
$$
\|Tx\|_{Y_0} \leqslant M_0\|x\|_{X_0}\ \ \forall \ x\in X_0\quad \text{and}\ \
\|Tx\|_{Y_1} \leqslant M_1\|x\|_{X_1}\ \ \forall \ x\in X_1
$$
it follows that for any $0 < \theta <1,$ $T$ maps $X_\theta$ into $Y_\theta$ and
$$
\|Tx\|_{Y_\theta} \le M_0^{1-\theta} M_1^\theta \|x\|_{X_\theta}\ \ \forall \ x\in X_\theta.
$$

\noindent An important problem is to identify classes of nonlinear maps for which (a version of)  the above interpolation property holds. It was investigated by many authors, see e.g.  the papers
\cite{CaZ}, \cite{Cw}, \cite{Ga}, \cite{Li},  \cite{Ma1}, \cite{Ma2}, \cite{MPW}, \cite{Pe} as well as the books \cite{Ma2}, \cite{Tr}
and references therein. However, all these results are obtained under the assumption that a nonlinear map is well defined
on the entire spaces $X_0$ and $X_1$. In applications, such an assumption is often too restrictive. The only exception,
to our  knowledge, is the paper of Tartar \cite[Theorem 1]{Ta}.
Using the real interpolation method (constructed with the help of the $K$-functional, see \eqref{K}) and a setup where $X_0 = X_+$ and $Y_0 = Y_+$,  Tartar  proved the
interpolation property for a  class of nonlinear maps $\Phi: U \to Y_0,$
defined on an open, nonempty subset $U \subset X_0$ with the following properties:

\vspace{0.1cm}

\noindent (T1)  $\Phi$ is locally
$\alpha-$H\"older continuous for some given $0 < \alpha \le 1$;

\vspace{0.1cm}

\noindent (T2) For any $x\in U$ there exist a neighborhood $V\subset U$ of $x$ and a constant $c > 0$ so that for any $y \in V \cap X_1,$
\[
\| \Phi(y) \|_{Y_1} \le c(\|y\|^{\beta}_{X_1} +1)
\]
 for some given $\beta > 0$.

\vspace{0.1cm}

\noindent   Note that the set $V \cap X_1$ might have an infinite diameter
in $ X_1$ even if $V$ is a (small) ball in $X_0.$
Hence in applications, property (T2) is often not satisfied or difficult to verify.
It is this fact that motivated our study on nonlinear interpolation in \cite{KT1}, \cite{SS1}, \cite{SS4}
where applications to inverse problems of spectral theory were considered.
Within the setup of the complex method of interpolation, in \cite{KT1},
 \cite{SS1},
we established the interpolation property for a class of nonlinear maps, defined on balls with center at the origin,
satisfying assumptions which in many applications are rather easy to prove. In this paper, our aim is to extend the results of \cite{KT1}, \cite{SS4} to a more general setup.

%%%%%%%%%%%%%%%%%%%%%%%%%%%%%%%%%%%%%%%%%%%%%%%%%%%%%%%%%%%%%%%%%%%%%%%%%%%%%%%%%
%%%%%%%%%%%%%%%%%%%%%%%%%%%%%%%%%%%%%%%%%%%%%%%%%%%%%%%%%%%%%%%%%%%%%%%%%%%%%%%%%

\section{Results}\label{Results}

\noindent To state our results, we first need to introduce some additional notation and establish some auxiliary results.
For any given complex Banach space $X$, denote by $B(r, X)$ the open ball in $X$ of radius $r$ centered at the origin.
Throughout this section, let $(X_0, X_1)$ be a complex Banach couple with $X_1\subseteq X_0$ and
$X_\theta =[X_0, X_1]_\theta$, $0< \theta <1,$ be the complex Banach spaces, constructed by the complex method of interpolation.
We have the following auxiliary lemma.
%%%%%%%%%%%%%%%%%%%%%%%%%%%%%%%%%%%%%%%%%%%%%%%%%%%%%%%%%%%%%%%%%%%%%%%%%%%%%%%%%
\begin{lemma}\hspace{-2mm}{\bf .}\label{lem:embeddings}
Assume that the norms of the Banach spaces $X_1\subseteq X_0$ satisfy
\begin{equation}\label{norm}
\|x\|_{X_0} \leqslant c\|x\|_{X_1} \quad \forall x\in X_1
\end{equation}
for some positive constant $c>0$. Then for any $0 < \theta < 1,$
\begin{equation}\label{eq:X_theta-embeddings}
\mbox{(i)} \,\, \|x\|_{X_0}\le c^\theta\|x\|_{X_\theta}\,\,\, \forall x\in X_\theta\, ;\quad \mbox{(ii)}\,\,
\|x\|_{X_\theta}\le c^{1-\theta}\|x\|_{X_1} \,\,\, \forall x\in X_1\,.
\end{equation}
In particular,
\begin{equation}\label{balls}
B(c^{-1} r, X_1)\subseteq B(c^{-\theta} r, X_\theta) \subseteq B(r, X_0).
\end{equation}
\end{lemma}
%%%%%%%%%%%%%%%%%%%%%%%%%%%%%%%%%%%%%%%%%%%%%%%%%%%%%%%%%%%%%%%%%%%%%%%%%%%%%%%%%%
\noindent{\em Proof of Lemma \ref{lem:embeddings}.}
Inequality (i) of \eqref{eq:X_theta-embeddings} follows from the interpolation property
of linear operators.  Indeed, by \eqref{norm} the identity operators
$$
I: X_1 \to X_0 \quad\text{and} \quad I: X_0\to X_0
$$
are bounded by $c$ and 1, respectively. Hence,
$$
I: X_\theta =[X_0, X_1]_\theta \to [X_0, X_0]_\theta =X_0
$$
is bounded by $c^\theta$.
To prove inequality (ii) of \eqref{eq:X_theta-embeddings}, consider for any given
$x\in X_1$ the analytic function $f : \mathbb{C}\to X_1$, $f(z):=c^{z-\theta} x$.
Clearly $f\in \mathcal{H}(X_0,X_1)$. Furthermore, in view of the definition of the norms
$\|\cdot\|_{X_\theta}$, $\|f\|_{\mathcal H}$ and of \eqref{norm},
\[
\qquad \quad \|x\|_{X_\theta}\le\|f\|_{\mathcal H}=
\max\{c^{-\theta}\|x\|_{X_0},c^{1-\theta}\|x\|_{X_1}\}\le c^{1-\theta}\|x\|_{X_1}. \qquad \qquad
\mbox{\finishproof}
\]
The first main result of this paper is the following
%%%%%%%%%%%%%%%%%%%%%%%%%%%%%%%%%%%%%%%%%%%%%%%%%%%%%%%%%%%%%%%%%%%%%%%%%%%%%%%%
\begin{Th}\hspace{-2mm}{\bf .}\label{th:main}
Assume that the norms of the Banach spaces $X_1\subseteq X_0$ satisfy
\begin{equation}\label{eq:norm1}
\|x\|_{X_0} \leqslant c\|x\|_{X_1}
\quad \forall x\in X_1
\end{equation}
 for some positive constant $c>0$.
Let $(Y_0, Y_1)$ be an arbitrary Banach couple and for some $r>0,$
\begin{equation}\label{est}
\Phi:  B(r, X_0) \to Y_+
\end{equation}
be an analytic map with
 $\Phi (B(r, X_0)) \subseteq Y_0$ and $\Phi (B(c^{-1}r, X_1)) \subseteq Y_1 $ so that
\begin{equation}\label{est0}
\Phi:  B(r, X_0) \to Y_0
\end{equation}
and
\begin{equation}\label{est1}
\Phi|_{B(c^{-1}r, X_1)} :  B(c^{-1}r, X_1) \to Y_1
\end{equation}
are continuous and bounded by the constants $M_0$ and $M_1$ respectively,
\begin{equation}\label{eq:boundedness2}
\sup_{x\in B(r,X_0)}\|\Phi(x)\|_{Y_0}\le M_0 \quad\text{and} \quad \sup_{x\in B(c^{-1}r,X_1)}\|\Phi(x)\|_{Y_1}\le M_1\,.
\end{equation}
Then for any $0 < \theta < 1$,  $B(c^{-\theta} r, X_\theta) \subseteq B(r, X_0)$ and
$\Phi ( B(c^{-\theta}r, X_\theta)) \subseteq Y_\theta$. Furthermore,
\begin{equation}\label{est2}
\Phi|_{B(c^{-\theta}r, X_\theta)} :  B(c^{-\theta}r, X_\theta) \to Y_\theta
\end{equation}
is bounded. More precisely, for any $x\in B(c^{-\theta}r, X_\theta)$,
\begin{equation}\label{estint}
\|\Phi(x)\|_{Y_\theta} \leqslant M_0^{1-\theta}{M_1}^\theta\,.
\end{equation}
\end{Th}
%%%%%%%%%%%%%%%%%%%%%%%%%%%%%%%%%%%%%%%%%%%%%%%%%%%%%%%%%%%%%%%%%%%%%%%%%%%%%%%%%%
\noindent{\em Proof of Theorem \ref{th:main}.}
By Lemma \ref{lem:embeddings}, for any $0 < \theta < 1,$ $B(c^{-\theta} r, X_\theta) \subseteq B(r, X_0)$. Hence the map $\Phi$ is well defined on $B(c^{-\theta}r, X_\theta)$.
It remains to prove that $\Phi(B(c^{-\theta}r, X_\theta))\subseteq Y_\theta$ and that \eqref{estint} holds.
Take an arbitrary $x\in B(c^{-\theta}r,X_\theta)$. By the definition of the norm $\|x\|_\theta$ there exists
a function $f\in \mathcal H(X_0, X_1)$ such that
\[
f(\theta) = x\quad \text{and}\quad\|f\|_{\mathcal H} < c^{-\theta}r\,.
\]
In particular, in view of the definition of the norm $\|f\|_{\mathcal H}$, one has for any $t\in\mathbb{R}$,
\begin{equation}\label{eq:f}
\|f(it)\|_{X_0}< c^{-\theta}r\quad\text{and}
\quad\|f(1+it)\|_{X_1} < c^{-\theta}r\,.
\end{equation}
Consider the function $g : {\overline S}\to X_+, z \mapsto g(z) := c^{\theta - z} f(z)$.
Clearly, $g\in\mathcal{H}(X_0,X_1)$ and $x=g(\theta)$.
As $X_0$ and $X_+$ coincide and the norms of $X_0$ and $X_+$ are equivalent,
$g : {\overline S}\to X_0$ is continuous and bounded, and $f|_{S} : S\to X_0$ is analytic.
Moreover, in view of  \eqref{eq:norm1} and \eqref{eq:f}, for any $t\in\mathbb{R}$,
\begin{equation}\label{second}
\|g(1+it)\|_{X_0}=c^{\theta -1}\|f(1+it)\|_{X_0}  \leqslant  c^\theta \|f(1+it)\|_{X_1} <r\,
\end{equation}
and
\begin{equation}\label{first}
\|g(it)\|_{X_0}\leqslant  c^\theta\|f(it)\|_{X_0} <r .
\end{equation}
Applying Hadamard's three line theorem \cite[Lemma 1.1.2]{BL} then yields
\begin{equation}\label{eq:r-ball}
\|g(z)\|_{X_0}< r \quad \forall z\in{\overline S}\,.
\end{equation}
In addition, we get from \eqref{eq:f} that
\begin{equation}\label{third}
\|g(1+it)\|_{X_1}=c^{\theta-1}\|f(1+it)\|_{X_1} <  c^{-1} r\,.
\end{equation}
Inequality \eqref{eq:r-ball} allows to define $F : {\overline S}\to Y_+$,
\[
F(z) := M_0^{z-1} M_1^{-z} \Phi(g(z)).
\]
As $g : {\overline S}\to B(r, X_0)\subseteq X_0$ is continuous and,
by assumption, $\Phi:  B(r, X_0) \to Y_0 $ is bounded and continuous
we conclude from
the continuity of the embedding $Y_0\subseteq Y_+$ that  $F : {\overline S}\to Y_+$ is continuous and bounded.
Similarly, as $g|_S : S\to B(r, X_0)\subseteq X_0$ and, by assumption, $\Phi:  B(r, X_0) \to Y_+$ are analytic, the function $F|_S : S \to Y_+$ is analytic as well.
The inequalities \eqref{first} and \eqref{third} together with the bounds (\ref{eq:boundedness2})
as well as the continuity and the boundedness
of the maps $\Phi:  B(r, X_0) \to Y_0 $ and
$\Phi|_{B(c^{-1}r, X_1)} :  B(c^{-1}r, X_1) \to Y_1$
%of \eqref{est0} and \eqref{est1}
 then imply that
$F\in \mathcal H(Y_0, Y_1)$ and $\|F\|_\mathcal H \leqslant 1$.
Furthermore, as $F(\theta) = M_0^{\theta -1}M_1^{-\theta} \Phi(x)$ one concludes that
\[
\|\Phi(x)\|_{Y_\theta} =  M_0^{1-\theta} M_1^\theta \|F(\theta)\|_ {Y_\theta}
\le M_0^{1-\theta} M_1^\theta \|F(z)\|_\mathcal H \leqslant  M_0^{1-\theta} M_1^\theta.
\]
This completes the proof of Theorem \ref{th:main}. \finishproof

If the nonlinear maps $\Phi:  B(r, X_0) \to Y_0 $ and $\Phi|_{B(c^{-1}r, X_1)} :  B(c^{-1}r, X_1) \to Y_1$  in Theorem \ref{th:main} admit polynomial bounds of order $n$ for some $n \ge 1$, then the following theorem asserts that
the same is true for the maps
$\Phi|_{B(c^{-\theta}r, X_\theta)} :  B(c^{-\theta}r, X_\theta) \to Y_\theta$ with arbitrary $0 < \theta < 1$.
Actually, the following theorem generalizes the corresponding result of \cite {SS4}.
%%%%%%%%%%%%%%%%%%%%%%%%%%%%%%%%%%%%%%%%%%%%%%%%%%%%%%%%%%%%%%%%%%%%%%%%%%%%%%%%
\begin{corollary}\hspace{-2mm}{\bf .}\label{coro:main}
Assume that all the assumptions of Theorem \ref{th:main} are satisfied except that the boundedness condition \eqref{eq:boundedness2}  is replaced
for some integer $n \ge 0$ by
\begin{equation}\label{eq:stronger1}
\|\Phi(x)\|_{Y_0} \leqslant  M_0 \|x\|_{X_0}^n \quad \forall x\in B(r, X_0), \quad
\end{equation}
\begin{equation}\label{eq:stronger2}
\|\Phi(x)\|_{Y_1} \leqslant  M_1 \|x\|_{X_1}^n \quad \forall x\in B(c^{-1}r, X_1).
\end{equation}
Then for any $0\le \theta\le 1$ and any $x\in B(c^{-\theta}r, X_\theta)$,
\begin{equation}\label{int.6}
\|\Phi(x)\|_{Y_\theta} \leqslant  M^{1-\theta}_0 M_1^\theta \|x\|^n_{X_\theta}.
\end{equation}
\end{corollary}

%%%%%%%%%%%%%%%%%%%%%%%%%%%%%%%%%%%%%%%%%%%%%%%%%%%%%%%%%%%%%%%%%%%%%%%%%%%%%%%
\noindent {\em Proof of Corollary \ref{coro:main}.}
Let $x\in B(c^{-\theta}r, X_\theta)$, and set $\rho:= \|x\|_{X_\theta}$.
Choose $\varepsilon >0$ so that $r_1: = c^{\theta}(\rho +\varepsilon) <r$.
Then by (\ref{eq:stronger1}),
\[
\|\Phi(y)\|_{Y_0}
\leqslant M_0r_1^n,\ \  \forall\  y\in B(r_1, X_0)
\]
and by (\ref{eq:stronger2}),
\[
\|\Phi(y)\|_{Y_1} \leqslant M_1 (c^{-1}r_1)^n \ \ \forall\ y\in B(c^{-1}r_1, X_1).
\]
Theorem \ref{th:main} then implies that for all $y\in B(c^{-\theta}r_1, X_\theta)$,
\[
\|\Phi(y)\|_{Y_\theta }\leqslant M_0^{1-\theta} M_1^\theta(r_1)^{n(1-\theta)} (c^{-1}r_1)^{n\theta}=
M_0^{1-\theta} M_1^\theta (\rho+\varepsilon)^n.
\]
In particular, the inequality above holds for $y=x$. As $\varepsilon>0$ can be chosen arbitrarily small
Corollary \ref{coro:main} follows.   \finishproof

Our second main result says that the map
$\Phi|_{B(c^{-\theta}r, X_\theta)} :  B(c^{-\theta}r, X_\theta) \to Y_\theta$
of \eqref{est2} in Theorem \ref{th:main}
%(and Corollary  \ref{coro:main})
is analytic.

\begin{Th}\hspace{-2mm}{\bf .}\label{th:analyticity}
The map
$\Phi|_{B(c^{-\theta}r, X_\theta)} :  B(c^{-\theta}r, X_\theta) \to Y_\theta
$ of Theorem \ref{th:main} is analytic.
\end{Th}
{\em Proof of Theorem \ref{th:analyticity}.}
Inspired by arguments used in \cite{SS4}, we prove the claimed statement by showing that
$\Phi|_{B(c^{-\theta}r, X_\theta)}$ is represented
by a series of analytic maps which converges absolutely and uniformly in $Y_\theta$
on any ball $B(c^{-\theta}\rho, X_\theta)$ with $0 < \rho < r$ -- see \cite[Theorem 2 in Appendix A]{PT}.
First note that $Y_0$ and $Y_1$ are continuously embedded in $Y_+$, so  by
%Lemma \ref{lem:embeddings} and
Lemma \ref{lem:analyticity} in  Appendix,
the maps $\Phi:  B(r, X_0) \to Y_0 $ and $\Phi|_{B(c^{-1}r, X_1)} :  B(c^{-1}r, X_1) \to Y_1$
of \eqref{est0} and \eqref{est1} are analytic.
Being analytic, the map $\Phi : B(r,X_0)\to Y_0$ is represented by its Taylor's series at $0$
with values in $Y_0$,
\begin{equation}\label{eq:taylor}
\Phi(h)=\Phi(0)+\sum_{n=1}^\infty\Phi_n(h),\quad h\in B(r,X_0).
\end{equation}
Here $\Phi(0)\in Y_1$ and $\Phi_n(h) = \frac{1}{n !} d_0^n \Phi (h, \cdots,h)$, $n\ge 1,$
with $d_0^n \Phi$ denoting the n'th derivative of $\Phi$ at $0$ (cf (A4) in Appendix).
We remark that for any $n \ge 1$, $\Phi_n(h)$ is a bounded homogeneous polynomial of degree $n$ in $h$ with values in $Y_0$ and hence analytic, and that the series in (\ref{eq:taylor}) converges absolutely and uniformly in $h$ on any ball $ B(\rho, X_0)$ with $0 < \rho < r$.
Moreover, in view of Cauchy's formula (cf. Appendix), for any $h\in X_0$ with $ h \ne 0,$
\begin{equation}\label{eq:Phi_n}
\Phi_n(h)=\frac{1}{2\pi i}\oint_{|z|=\rho}\frac{\Phi(z h)}{z^{n+1}}\,dz
\end{equation}
where $\rho$ is chosen arbitrarily so that $0<\rho<r/\|h\|_{X_0}$.
Using \eqref{eq:Phi_n} and the first inequality in \eqref{eq:boundedness2} it follows that for any
$h\in X_0$ with $ h \ne 0,$
\begin{eqnarray*}
\|\Phi_n(h)\|_{Y_0}\le \frac{1}{2\pi}\int_0^{2\pi}\frac{\|\Phi(\rho e^{i t} h)\|_{X_0}}{|\rho e^{i t}|^{n+1}} |\rho\,i e^{i t}|\,dt
\le \frac{M_0}{\rho^n}\,.
\end{eqnarray*}
As this inequality holds for any $0<\rho<r/\|h\|_{X_0}$ one has for any $h\in X_0$ with $ h \ne 0,$
\begin{eqnarray}\label{eq:0-end}
\|\Phi_n(h)\|_{Y_0}\le\frac{M_0}{r^n}\,\|h\|_{X_0}^n\,.
\end{eqnarray}
Note that the latter estimate holds trivially for $h = 0$.
Applying the arguments above to the analytic map $\Phi|_{B(c^{-1}r,X_1)} : B(c^{-1}r,X_1)\to Y_1$
we see that for any $h\in B(c^{-1}r,X_1)$,  $\Phi_n(h) \in Y_1$ and
the series \eqref{eq:taylor} converges in $Y_1$. Furthermore by (\ref{eq:Phi_n}), for any $h\in X_1$,
\begin{eqnarray}\label{eq:1-end}
\|\Phi_n(h)\|_{Y_1}\le\frac{M_1}{(c^{-1}r)^n}\,\|h\|_{X_1}^n\,.
\end{eqnarray}
Applying Corollary \ref{coro:main} to the the map $\Phi_n : X_0\to Y_0\subseteq Y_+$ we obtain from
\eqref{eq:0-end} and \eqref{eq:1-end} that for any $h\in X_\theta$, $\Phi_n(h) \in Y_\theta$ and
\begin{equation*}
\|\Phi_n(h)\|_{Y_\theta}\le \left(\frac{M_0}{r^n}\right)^{1-\theta}\left(\frac{M_1}{(c^{-1}r)^n}\right)^\theta\|h\|_{X_\theta}^n
\le M_0^{1-\theta} M_1^\theta\left(\frac{\|h\|_{X_\theta}}{c^{-\theta}r}\right)^n\,.
\end{equation*}
This inequality shows that for any $h \in B(c^{-\theta}r,X_\theta)$, the series in \eqref{eq:taylor} converges
in $Y_\theta$ and that it converges
absolutely and uniformly on any ball $B(c^{-\theta} \rho,X_\theta)$ with $0 < \rho < r$.
\finishproof

\noindent Finally, we remark that Theorem \ref{th:main}, combined with Theorem \ref{th:analyticity},
generalizes Theorem 1.1 of \cite{KT1} in the context of the setup, chosen in this paper.

%%%%%%%%%%%%%%%%%%%%%%%%%%%%%%%%%%%%%%%%%%%%%%%%%%%%%%%%%%%%%%%%%%%%%%%%%%%%%%%%
%%%%%%%%%%%%%%%%%%%%%%%%%%%%%%%%%%%%%%%%%%%%%%%%%%%%%%%%%%%%%%%%%%%%%%%%%%%%%%%%

\section{Appendix}\label{Appendix}
In this appendix we review the notion of an analytic map between complex Banach spaces
and discuss properties of such maps needed in Section \ref{Results}.
For more details we refer the reader e.g. to \cite{Die} or \cite[Appendix A]{PT}.

\noindent Let $X$ and $Y$ be complex Banach spaces and let $U\subseteq X$ be an open set in $X$.
A map  $F : U \to Y$ is called {\em analytic} if it is Fr\'echet differentiable over $\mathbb C$ at any point
$x\in U$. The map $F : U \to Y$ is called {\it weakly analytic} if  for any $x\in U$, $h\in X$, and for any
 $f\in Y^*$, the complex-valued  function $f(F(x +zh))$ is holomorphic in a small disk  in $\mathbb C$
centered at zero.
% In the case where $X$ is a separable Hilbert space with orthonormal basis $\{e_n\}_{n\ge 1}$ and scalar product
%$\left<\cdot,\cdot\right>$, the weak analyticity is equivalent to the Gateaux differentiability over $\mathbb C$ of the %coordinate
%functions $\left<F(x), e_n\right>$.

\noindent Let us recall the following analyticity criteria (see e.g. \cite[Appendix A, Theorem 1.1]{PT}).
Assume that $X$ and $Y$ are complex Banach spaces and that $\Phi: U\to Y$ is a map defined on an open subset  $U$ of $X$, $U\subseteq X$.
Then the following statements are equivalent:

\vspace{0.1cm}

\noindent (A1) {\em $\Phi: U\to Y $ is analytic;}

\vspace{0.1cm}

\noindent (A2) {\em $\Phi: U\to Y $  is weakly analytic and locally bounded; }

\vspace{0.1cm}

\noindent (A3) {\em $\Phi: U\to Y $  is continuous and for any
$x\in U, h\in X, $ there exists a disk $D_r =\{ z\in \mathbb C \vert \ \, |z|<r\}$ such that
for any $0 < \rho < r,$
Cauchy's formula holds
\begin{equation}\label{integ}
F(z)=\frac1{2\pi
i}\int\limits_{|\xi|=\rho}\!\!\frac{F(\xi)}{\xi-z}\,\,d\xi,\quad\forall \, |z|< \rho;
\end{equation}

\vspace{0.1cm}

\noindent (A4) {\em $\Phi: U\to Y $  is infinitely differentiable on $U$ and is represented by
its Taylor series in a neighborhood of each point of $U$, i.e., for any $x \in U$,
$\Phi(x+h)=\Phi(x)+\sum_{n=1}^\infty d_x^n\Phi (h, \cdots ,h)$ where
$d_x^n\Phi$ denotes the $n'th$ derivative of $\Phi$ at $x$ and the
series converges absolutely and uniformly for any $h$ in $B(\rho, X)$ with  $0 < \rho < r$
so that $B(r, X) \subset U$.}

\vspace{0.1cm}
}

%%%%%%%%%%%%%%%%%%%%%%%%%%%%%%%%%%%%%%%%%%%%%%%%%%%%%%%%%%%%%%%%%%%%%%%%%%%%%%%%%%%
\noindent Now let us prove the following
\begin{lemma}\label{lem:analyticity}
Let $X$, $Y$ and $Y_+$ be Banach spaces so that $Y\subseteq Y_+$ is continuous.
Let $U\subseteq X$ be an open set in $X$ and let $\Phi : U\to Y_+$ be analytic, $\Phi(U)\subseteq Y$,
and $\Phi : U\to Y$ be continuous. Then $\Phi : U\to Y$ is analytic.
\end{lemma}

\noindent{\em Proof of Lemma \ref{lem:analyticity}.}
For any $x\in U$, $h\in X$,  consider the map $F: D_r \to Y_+$, $F(z):=\Phi(x+zh)$,
where the  disk $D_r\subseteq\mathbb C$ is centered at zero and chosen so that $x+z h \in U_r$ for any $z\in D$. As by assumption, $F$ is analytic in the $Y_+$-norm
formula \eqref{integ} holds in $Y_+$ for any given $z\in D$.
Furthermore, as we assume in addition that $F$ is continuous with respect to the $Y$-norm,
the integral on the right hand side of \eqref{integ}
defines a continuous map $G: D \to Y$. In view of the embedding $Y\subseteq Y_+$, one has $G=F$.
\finishproof

\bigskip

\noindent T. Kappeler, Institute of Mathematics, University of Z\"urich, Z\"urich,\\
Switzerland. Email: thomas.kappeler@math.uzh.ch

\hspace{1cm}

\noindent A. Savchuk, Department of Mechanics and Mathematics, Lomonosov Moscow\\
 State University, Moscow, Russia. Email: artem\_savchuk@mail.ru

\hspace{1cm}

\noindent A. Shkalikov, Department of Mechanics and Mathematics, Lomonosov Moscow
State University, Moscow, Russia. Email: ashkalikov@yahoo.com

\hspace{1cm}

\noindent P. Topalov, Department of Mathematics, Northeastern University,\\
Boston, Massachusetts, USA. Email: p.topalov@neu.edu


\begin{thebibliography}{99}

\bibitem{BL} J. Bergh, J. L\" ofstr\"om, {\it Interpolation spaces. An introduction}, Springer, 1976

\bibitem{CaZ} A. Calder\'on, A. Zygmund, {\it  A note on the interpolation of sublinear operations},
Amer. J. Math., $\bf 78$(1956),  282-288

\bibitem{Cw} M. Cwikel, {\it A counterxample in nonlinear interpolation}, Proc. Amer. Math. Soc., $\bf 62$(1977), 62-66

\bibitem{Die} J. Dieudonne, {\it Foundations of Modern Analysis}, 5th ed., Academic Press, 1960

\bibitem{Ga} E. Gagliardo, {\it Interpolation d'espaces de Banach et applications III},
 C.R. Acad. Sci. Paris, $\bf 248$(1959), 3517-3518

\bibitem{KT1} T. Kappeler, P. Topalov, {\it  On nonlinear interpolation}, to appear in Proc. Amer. Math. Soc.,
ArXiv:1306.5721

\bibitem{Li} J. Lions, {\it Interpolation lin\'eaire et non lin\'eaire et regularit\'e,}
Instituto Nazionale di Alta Matematica, Symposia Math., $\bf 7$(1971), 443-458

\bibitem{Ma1} L. Maligranda, {\it On interpolation of nonlinear operators},
Annales Societatis Mathematicae Polonae. Series 1: Commentationes Mathematicae, $\bf 28$(1989),  253-275

\bibitem{Ma2} L. Maligranda, {\it A Bibliography on "Interpolation of Operators and Applications": (1926-1990)},
H\"ogskolan i Lulea, Lulea Univ., 1990.

\bibitem {MPW} L. Maligranda, L. Persson and  J. Wyller,  {\it Interpolation and partial differential
equations}, J. Math. Phys., $\bf 35$(1994), no. 9, 5035-5046

\bibitem {Pe}  J. Peetre, {\it Interpolation of Lipschitz operators and metric spaces},  Mathematica (Cluj),
$\bf 12$(1970),  325-334

\bibitem{PT} J. P\"oschel, E. Trubowitz,  {\it Inverse Spectral Theory},  Academic Press, 1987

\bibitem{RS} M. Reed, B. Simon, {\em Methods of Modern Mathematical Physics II},
Academic Press, 1975

\bibitem {SS1} A. Savchuk, A. Shkalikov, {\it On the eigenvalues of the Sturm-Liouville
operator with potentials from Sobolev spaces}, Mathematical Notes,  $\bf 80$(2006), no 6, 814-832

\bibitem {SS2}  A. Savchuk, A. Shkalikov, {\it On the properties of maps associated with inverse Sturm-Liouville problems},
Proceedings of Steklov Math. Institute, $\bf 260$(2008), 218-237

\bibitem{SS3} A. Savchuk, A. Shkalikov, {\it Inverse Sturm-Liouville problems with potential in Sobolev spaces. Uniform stability},
Funct. Anal and its Apppl., $\bf 43$(2010), no 3, 270-285

\bibitem{SS4} A. Savchuk, A. Shkalikov, {\em On the interpolation of analytic mappings}, Mathematical Notes, $\bf 94$(2013), no 4,
547-550,  ArXiv:1307.0623.

\bibitem{Ta} L. Tartar {\it Interpolation non lin\'eaire et r\'egularit\'e},
J. of Funct. Analysis, $\bf 9$(1972), 469-489

\bibitem{Tr} H. Triebel, {\it Interpolation Theory, Function Spaces, Differential Operators},
 North-Holland, 1978

\end{thebibliography}
\end{document}